\documentclass[11pt,twoside,a4paper,notitlepage]{article}

\usepackage{graphicx}
\usepackage{mathtools}
\usepackage{enumerate}
\usepackage{amssymb}
\usepackage{amsthm}
\usepackage{mathptmx} 
\usepackage{csquotes}
\usepackage{latexsym}
\usepackage{tikz}
\usepackage[english]{babel}
\usepackage[T1]{fontenc} 
\usepackage{enumerate}

\theoremstyle{definition}

\newtheorem{definition}{Definition}

\newcommand{\impl}{\Rightarrow}

\begin{document}

\title{Infinity and Continuum in the Alternative Set Theory}

\author{Kate\v{r}ina Trlifajov\'{a}}

\title{Sizes of Countable Sets}

\author{Kate\v{r}ina Trlifajov\'{a}}

%\affiliation{Czech Technical University in Prague, Faculty of Information Technology}

%\email{katerina.trlifajova@fit.cvut.cz}

 \date{}

\maketitle

\maketitle

\begin{abstract}

Alternative set theory was created by the Czech mathematician
Petr Vop\v enka in 1979 as an alternative to Cantor's set theory. Vop\v enka criticised
Cantor's approach for its loss of correspondence with the real
world. Alternative set theory can be partially axiomatised and regarded
as a nonstandard theory of natural numbers. However, its intention is much wider. It attempts to retain a correspondence between mathematical notions
and phenomena of the natural world. Through infinity, Vop\v enka grasps
the phenomena of vagueness. Infinite sets are defined
as sets containing proper semisets, i.e. vague parts of sets limited by the horizon. The new interpretation extends the field of applicability of mathematics and simultaneously indicates its limits. 
This incidentally provides a natural solution to some classic philosophical problems such as the composition of a continuum, Zeno's paradoxes and sorites.
Compared to strict finitism and other attempts at a reduction of the infinite to the finite Vop\v enka's theory reverses the process: he models the finite in the infinite. 


\end{abstract}

\section{Introduction }\label{intro}

Vop\v enka's Alternative set theory (AST) represents an attempt to present a new set
theory on a phenomenal basis. It is \enquote{a mathematical-philosophical attempt in
the study of the infinite}. Both components are important and intertwined.

This is not the only reason why the explanation of Vop\v enka's theory is difficult. Its
mathematical form has developed and changed several times. Vop\v enka
repeatedly searched for a better way to grasp fundamental phenomena
that he thought rational science had overlooked. He gradually developed his
theory, left blind alleys and looked for new ones.

Although common set-theoretical notions are used, the framework of AST is not Zermelo-Fraenkel set theory (ZF), and Vop\v enka introduces new, unusual concepts that need to be justified. 

In the 1960s, Vop\v enka had collaborated with an international
mathematical community. He had to interrupt it after 1968
when Czechoslovakia was closed behind the Iron Curtain. However, he continued his
research with a small group of his colleagues in the 1970s and 1980s. After the Velvet Revolution in 1989, he became involved in politics.\footnote{He was the vice-rector of Charles University in Prague in 1990 and the minister of education in 1990 -- 1992.} Vop\v enka returned to his research at the end of the millennium. Unfortunately, he had not re-established his relationship with the
international community. It is a pity since it could have been mutually enriching.
Much of his work has never been translated and 
exists only in Czech.

The original version of Vop\v enka's AST is nearly fifty years old. The question
is whether it is still worthwhile to consider it from anything other than a historical
point of view.

AST tries to bridge a gap between
infinite mathematical objects and finite physical entities. Its philosophical justification 
is distinctive and offers a new perspective that may inspire
further research. The new interpretation extends the field of applicability of mathematics while suggesting its limits. The themes that Vop\v enka dealt with reappear with new intensity at
the beginning of the new millennium, particularly in connection with vagueness theory, feasibility and alternatives to the classical set theory (Fletcher 2007; Bellotti 2008; Gaifman 2010; Dean 2018; Holmes, Forster, Libert 2012).

\subsection{Overview}

The main aim of this paper is to present a phenomenological conception of the continuum
and the infinite in AST and its evaluation in the present context.

For a short history of AST and its context, see Section \ref{hist}. Vop\v enka justified the creation
of his theory through the critique of Cantor's set theory, \ref{Cantor}. Some of the predecessors
and contemporaries that influenced Vop\v enka are mentioned in \ref{related}.
The development of Vop\v enka's theories is briefly described in \ref{evolution}.

The basic mathematical notions of AST and their philosophical justification are presented in the Sections
\ref{inf} - \ref{cont}. An explanation of the concept
of the infinite appears in \ref{inf}, and one for the continuum is provided in \ref{cont}.
Number structures are described in \ref{numb}. Equipped with new concepts, AST allows to describe some phenomena that classical mathematics cannot capture, see \ref{sigma}, \ref{indiscernibility}, \ref{time} and \ref{motion}. Natural solutions to some classical Greek paradoxes are mentioned throughout the text. 

The last two sections contain my reflexion of AST. Section \ref{ideal} concerns a discussion on abstraction and idealization in AST. Section \ref{concl} is a conclusion.

\section{Historical context}\label{hist}

\subsection{Critique of Cantor Set Theory}\label{Cantor}
Vop\v enka considered his set theory as an alternative to Cantor's set theory, which he
criticised for several reasons. Since Cantor's set theory assumes the existence of
actually infinite sets, it loses touch with the real world and becomes motivated
only by formal considerations.

The existence of undecidable theorems, such as the Axiom of choice or the
Continuum hypotheses divide set theory into several branches,
none of which can be viewed as being true.

On one side, set theory opened the way for researchers to study an immense number
of various infinite structures. On the other side, it closed the door to the study
of vague structures that mathematics is unable to grasp.

Later, Vop\v enka also pointed out Cantor's theological reasoning regarding actual
infinity (Dauben 1990, pp. 228 - 232). Mathematicians took a standpoint that can be described
as the \enquote{God-like position} by claiming the right to decide what is true
about infinite sets.

\subsection{Related Topics}\label{related}
It is no accident, and perhaps it was the atmosphere of the time in which AST was created that similar topics had been investigated.\footnote{Zuzana Hanikov\' a deals in detail with the contexts and influences of other theories on AST and its development. (Hanikov\' a 2022).} It is far beyond
the scope of this paper to review them all, I will just mention a few that influenced Vop\v enka. 

\subsubsection{Non-standard Analysis}

From the formal and technical point of view, AST can be considered as a particular case of \emph{Nonstandard Analysis} (NSA). (Robinson 1966). NSA is founded on Cantor set theory that enriches with new techniques, whereas AST is partially an informal theory based on a different philosophical background.

In this time, Nelson's \emph{Internal Set Theory} (Nelson 1977) was created. That is an axiomatic/syntactic variant of NSA without any special philosophical reasoning. Vop\v enka's former student, Karel Hrb\' a\v cek, axiomatically developed his own
\emph{Non-Standard Set Theory} that is also relatively consistent with ZF. (Hrb\' a\v cek, 1979).
However, this was not the way Vop\v enka wanted to go. 

\subsubsection{Feasible numbers}\label{feas}
Another source of inspiration was Paul Bernays's paper \emph{On Platonism in Mathematics} where Bernays noted a discrepancy between the representation of
numbers in the Arabic system and its concrete accessibility, for instance, of the number $67^{257^{729}}$. (Bernays 1936, pp. 12 - 13).

He suggested \emph{strict finitism} as a conceivable position. Strict finitists are
concerned with the assumption that some natural numbers cannot be achieved
by simple arithmetic operations like a successor function. 

This was followed by Alexander
Yessenin-Volpin who declared his adherence to strict finitism and sketched a
program for proving the consistency of Zermelo-Fraenkel set theory in strict
finitist mathematics in the late 1950s.\footnote{\enquote{How to formalise the intuitive notion of feasible numbers? To see what feasible numbers are, let us start by counting: $0,1,2,3,$ and so on. At this point, Yessenin-Volpin (in his \emph{analysis of potential feasibility}, 1959) asks: \enquote{What does this 'and so on' mean? Up to what extent 'and so on'?} And he answers: \enquote{Up to exhaustion!} Note that by cosmological constraints exhaustion must occur somewhat before, say, $2^{1000}$.} (Sazonov 1995, p. 30).}

Rohit Parikh also 
dealt with the question of \emph{feasibility}
of concrete functions and procedures. 
He demonstrated that formal systems in which
\enquote{large} numbers are treated as if they were infinite provide correct results for all
proofs of a \enquote{reasonable} length. Parikh investigated the possibility of replacing infinite
numbers with large finite numbers denoted by primitive recursive terms.
He added a new predicate $F$, \emph{feasible}, to Peano arithmetic such that $F(0)$ and $F$ was closed on arithmetic operations. A number $2_{2^{1000}}$ constructed recursively via repetitive exponentiation was not feasible. 
This theory
is \enquote{almost consistent} in the sense that every formal proof of a contradiction
should contain at least $2^{1000}$ symbols. (Parikh 1971).

Parikh's ideas were later developed by Vladimir Sazonov. He considered Parikh's \enquote{large number} 
too rough upper bound for feasible numbers
and proposed another system, FEAS, to formalise vague notions such as feasible
numbers. Sazonov's restriction on proofs consists of allowing only normal,
natural deductions\footnote{In particular, this means that one cannot freely use the general modus ponens rule $E$ with the corresponding rule $I$ in the system of Natural deduction calculus. Then the implication does not have to be transitive.}. Moreover, the number of symbols in proof should be
(intuitively) feasible. This system is consistent and entails that $2^{2^{10}} = \infty$. (Sazonov 1995).

The problem in question was the existence of non-trivial cuts on natural numbers. Is there a subset (a subclass) of natural numbers closed on predecessors that does not have the greatest element and is upper bounded by a concrete number? These should be feasible numbers. 

\emph{Wang's paradox} represents the same problem: \enquote{The number $0$ is small; if $n$ is small then $n + 1$ is small. Therefore, every number is small.} There are many such examples, variants of the ancient paradox of a heap. After a careful analysis, Michael Dummett concluded that they are always connected with a vague predicate and demonstrated that such predicates
are intrinsically inconsistent. So, strict finitism is an untenable position. (Dummett 1975)

Dummett's repudiation was later criticised. Dean argued that it is
based on assumptions about number systems that strict finitists would almost certainly have rejected.
(Dean 2018, p. 295). Dummett's paper is
otherwise interesting due to his thorough analysis of vagueness.

\subsubsection{Analysis without Actual Infinity}\label{Mycielski} 
Supporting Vop\v enka was the paper of Jan Mycielski \emph{Analysis without Actual Infinity}. (Mycielski 1988).  Mycielski defined a simple axiomatisation FIN of a first-order theory such
that every finite part of FIN has finite models. FIN is sufficient for the
development of the analysis in the respect that all applicable mathematical theorems
can be correctly translated and proved in
FIN. However, Mycielski warns that the formalisation of analysis in FIN is clumsy.

\subsection{Evolution of Vop\v enka 's theories}\label{evolution}

\subsubsection{Vop\v enka's Principle}
Petr Vop\v enka (1935 -- 2015) achieved significant results as a set theorist in the
1960s. He became well-known for his invention of a large cardinal axiom, which is now
called Vop\v enka 's principle, in 1965. This principle is still used (Gitman,
Hamkins, 2018). It is stronger than the existence
of measurable cardinals and has several equivalent formulations, one of
which is that every proper class of first-order structures contains two different
members. One of them can be elementarily embedded in the other. See
(\v Svejdar, p. 1266).\footnote{According to (Pudl\' ak 2013, p. 204), Vop\v enka 's principle was originally intended as a joke: Vop\v enka was apparently unenthusiastic about large cardinals and introduced his principle as a bogus large cardinal property, planning to show later that it was not consistent. However, before publishing his inconsistency proof, he found a flaw in it.}

\subsubsection{Theory of Semisets}
In 1972, Vop\v enka and H\' ajek, published the book \emph{Theory of Semisets} (H\' ajek \& Vop\v enka 1972), TS. They had already proved the existence of non-trivial definable cuts in G\H odel-Bernays set theory (GB). (H\' ajek \& Vop\v enka 1973).
Roughly speaking, one can think of a model of GB as a model of ZF where some subsets are interpreted as classes.
These \enquote{subcollections} were called semisets. (\v Svejdar 2018).

TS had not been accepted with great enthusiasm. It was written in a highly formal language, and its results were rather incomprehensible.\footnote{Azriel L\' evy wrote in his review: \enquote{It was far enough to convince readers that modern metamathematics can be carried out for TS to the extent that it is carried out for ZF. 
As a result of the unusual way TS handles set theory and, even more, as a result of the highly formal approach taken in writing this book, the wealth of information in it is almost completely inaccessible to the students of set theory. This is a pity since the book contains many of the results of Vop\v enka 's Czech school of set theory and shows how to obtain the independence proofs of set theory by means of relative interpretations.} (L\' evy, p. 1423). } 

The theory had not been further extended in this form. However, both H\' ajek and Vop\v enka continued to research vagueness, but each went in a different direction. Later, H\' ajek published a book \emph{Metamathematics of Fuzzy Logic}, \enquote{a systematic treatment of deductive aspects and structures of fuzzy logic understood as many-valued logic \textit{sui generis}}. (H\' ajek 1998).

\subsubsection{Development of the Alternative Set Theory}

The main principles of AST had already been formulated in 1974. After that, Vop\v enka further developed it in a seminar that he led. The book \emph{Mathematics in the AST} (Vop\v enka 1979) was published five years later. 
This book was the only complete English version.\footnote{ \label{AST} The original book (Vop\v enka 1979) is hardly available. A copy can be found at https://drive.google.com/file/d/17JRj2orUVDw7lrBEmBS1K6OK06RP32Xa/view. (Holmes 2017) gives an abbreviated overview of its axioms.}

Except for the introduction, \emph{Mathematics in the Alternative Set Theory} is mostly a mathematical
book. Explanations of new concepts and ideas are only briefly indicated.

From a purely formal point of view, an important part of AST can be
axiomatised. The resulting axiomatic system can be identified with a particular rather weak
version of nonstandard set theory.\footnote{Nonstandard models of Peano arithmetic contains infinite numbers. AST can be formally described as an $\omega$-saturated model of cardinality $\aleph_1$ of Peano arithmetics. Robinson's Nonstandard Analysis uses ultrafilters (Robinson 1966) to construct a model of nonstandard real numbers where the differential and integral calculus can be consistently described using infinitely small quantities.} This proves its relative consistency with respect to ZF.

Vop\v enka's collaborators A. Sochor, J. Ml\v cek, K. \v Cuda, A. Vencovsk\' a and others developed a formal mathematical perspective of AST. They investigated its consistency, its models in ZF, its extensions and its metamathematics.\footnote{\enquote{A model of AST can easily be constructed as follows. Let $HF$ be the set of hereditarily finite sets. Let $(\widetilde V, \widetilde E)$ be the ultrapower of $(HF, \in)$ over some nontrivial ultrafilter on $\omega$. Add to $(\widetilde V, \widetilde E)$ all subsets $X \subseteq \widetilde V$ such that for no $x \in \widetilde V$ do we have $X = \{y; (\widetilde V, \widetilde E) \models y \in x \}$. If we assume the continuum hypothesis, then the resulting model is a model of AST.} (Pudl\' ak, Sochor 1984, p. 572).}

But that was not exactly the direction Vop\v enka wanted to go. He considered AST a non-formalised \enquote{naive} theory like Cantor set theory. Although some of its important fragments can be axiomatised, it is more important to retain a correspondence between mathematical notions and natural phenomena. 

The extended Slovak version \emph{Introduction to Mathematics in the Alternative Set Theory} followed (Vop\v enka 1989). In this book, Vop\v enka developed both mathematical and philosophical aspects of AST. He explained the phenomenological justification for his new mathematical concepts: objects, classes, sets, semisets, equality, $\pi$-classes and $\sigma$-classes.

Vop\v enka clarified his philosophical approach in \emph{Meditations on
Foundations of Science} in detail (Vop\v enka 2001). He explained his newly introduced
concepts that describe the phenomena of the natural world that cannot be
grasped in classical mathematics, such as horizon, vagueness, indiscernibility.

\subsubsection{New Infinitary Mathematics}
Vop\v enka returned once again to the philosophical and mathematical
foundations of his theory in the new millennium. He tried to improve it so that it was
closer to the real world. The last version he completed was
published under the title \emph{New Infinitary Mathematics} (NIM) after his death in 2015.

While AST is partially axiomatisable with some
philosophical comments, NIM is even looser. Some principles are stated here so that they may be
interpreted both in nonstandard analysis and in applied
mathematics. Vop\v enka again returns to the key question of non-trivial cuts on natural numbers. NIM also contains advanced mathematical theorems, among them a problematic proof of the non-existence of the set of natural numbers 
using an ultra-power over a covering structure. 

NIM's programme is characterised by the following imperative:
\begin{quote} 
Wherever a vagueness occurs, look for a horizon and natural infinity that has caused this vagueness; then idealise this situation! (Vop\v enka 2015, p. 13). \end{quote}

NIM is more complicated and less comprehensible than AST. Vop\v enka wished to avoid
being committed to any binding axiomatic framework, and he even left space to interpret what logic should be used.

\subsubsection{The present version}
The main goal of this paper is to explain the philosophical reasoning of mathematical notions of AST (if we can call
philosophical reasons the correspondence between the real and mathematical
world). 
I will proceed from the
only complete English version \emph{Mathematics in the Alternative Set
Theory} (Vop\v enka 1979) and simply describe several important concepts, their
justifications and relations. 
I consider the phenomenological explanations in other works, especially in (Vop\v{e}nka 2001) and (Vop\v{e}nka 2015).

Of course, the present interpretation is influenced by my comprehension. I select the parts of AST that I consider revealing and valuable for further research. 
Some questions regarding Vop\v enka's theory have remained unanswered, and I will try to 
offer an answer.

\section{Infinity in AST}\label{inf}

\subsection{Vagueness}

Vop\v enka refers to Husserl's challenge to go \enquote{back to the things themselves}.
Doing this means returning to the ways that things are actually given in our experience.\footnote{Husserl presents his \emph{Principle of All Principles}: 
\enquote{Every originary presentive intuition is a legitimising source of cognition, that everything originally (so to speak, in its `personal' actuality) offered to us in `intuition' is to be accepted simply as what it is presented as being, but also only within the limits in which it is presented there} (Husserl 1982, sec. 24).}
He considered vagueness to be a fundamental phenomenon that rational
science has overlooked. The character of mathematics from the essence embedded in it in antiquity rules out vagueness from the mathematical world.
\begin{quote} In this world, there are no distinct boundaries of greenness, slowness, weakness, love, and so forth. (\dots)
Vagueness is a phenomenon that modern science attempts in vain to cast out from the maze of phenomena in which we are thrown. In it is the key to grasping the separation of the world of science from the natural world and to comprehending the incomprehensible phenomena of the natural world through modern science. (Vop\v{e}nka 2001, p. 56). \end{quote} 

\subsection{Horizon}

Every observation we make, even in any direction, is limited. 
Either it is obstructed by a fixed boundary, which sharply interrupts it, or it is limited by a horizon towards which the distinction of our view disappears. This applies not only to optical observations; the horizon is understood in a much broader figurative sense.

Vop\v enka borrowed the concept of a horizon 
from Husserl. He used it in a slightly different technical way, but for the same reasons. 
The horizon is a limit that separates the field of a direct experience from
that of an indirect experience. It is not a sharp boundary; the illuminated part passes continuously into the unilluminated one. The closer one gets to the horizon; the more vagueness surrounds the illuminated part. 

The horizon always depends on the observer and their point of view. Only
through our ability to reflect do we know of its existence in the natural world.\footnote{Vop\v enka's concept of a horizon resembles that one: \enquote{The horizon is tied to an observer; there is something subjective about it. On the other hand, it appears outside of man; it is eternal; we can never catch it (p. 183). 
The horizon adds nothing to the world. The horizon does not enrich the world. On the other hand, the world without the horizon is unimaginable, even impossible (p. 184). 
Man cannot remove it nor reach it. To man, the horizon represents the idea that there is more than what he sees (p. 187). (Van Peursen 1977).}} The horizon limits our view into the distance, but also into the depths, towards ever smaller objects.
\subsection{Classes, Sets and Semisets}

To describe such phenomena mathematically, Vop\v enka introduced the new notion of the semiset in
addition to the classical notions of set and class. 

\begin{itemize}
\item A \emph{class} is any collection of objects which we consider autonomous. 
\item A \emph{set} is a class that is distinctly defined. Its borders are precisely determined. Its elements could (hypothetically) be arranged in a list.
\item A \emph{semiset} is a vague part of a set or a subclass of a set that is not a set.
\end{itemize}

Proper classes are classes that are not set. They represent objects that
are bordered somewhere by a horizon. The vagueness of a proper class is present for a multitude of its 
elements. It is not a case of semantic indeterminacy when vagueness occurs due to
mere imprecision, ambiguity or obscurity of an object.

Semisets can be found nearly everywhere. Many examples of semisets are given in the literature concerning vagueness
and feasibility, Yessenin-Volpin's \enquote{feasible numbers},
Wang's \enquote{small numbers}, \enquote{orange objects} or \enquote{heartbeats in someone's childhood} (Dummett 1975), \enquote{inhabitants of a small town} or {age of a childish in hours} (Dean 2018), \enquote{walking distance} (Gaifman 2010).\footnote{Walking distance is still a walking distance if we increment it by one foot (but not 5 miles); a child is still a child 1 hour later (but not 5 years). (Gaifman 2010, p. 6).} 

Ancient \emph{sorites paradoxes} always contain a semiset. Let us consider the \enquote{bald man paradox}. Hairs one can pull out of a hairy man so that he is not bald is a vague part of the set of all his hairs. It is a semiset. 

\subsection{Axioms}
These definitions are rather informal. However, they are specified axiomatically.
AST uses common notions and symbols from the language of sets such as membership relation, equality, empty set, subclass, intersection, union, linear order, function, one-to-one function, isomorphism, etc.

A \emph{set formula} is a formula containing only set variables or constants, and all quantifiers are restricted to sets.
\medskip

\centerline{Axioms for sets}
\begin{enumerate}
\item[(A1)] \emph{Empty set:} There is an empty set $\emptyset$.
\item[(A2)] \emph{Successor:} If $x$ and $y$ are sets then $x \cup \{y\}$ is a set.
\item[(A3)] \emph{Extensionality:} Sets with the same elements are equal. 

\item[(A4)] \emph{Induction:} Let $\varphi$ be a set formula. If $\varphi(\emptyset)$ and 
$(\forall x)(\forall y)(\varphi(x) \impl \varphi(x \cup \{y\})$ then $(\forall x)\varphi(x)$. 
\item[(A5)] \emph{Regularity:} Every set has an element disjoint from it. 
\end{enumerate}

Formally, one can work in the \emph{universe of sets} that is constructed iteratively from the empty set according to the axioms. The universe of sets corresponds to hereditary finite sets in Cantor set theory.\footnote{For this reason, Vop\v enka claimed that \enquote{all sets are finite in Cantor's sense.} 
Sochor even asserts \enquote{There is no infinite set.} (Sochor 1984, p. 172). The reason is that he uses slightly different terminology and considers all sets to be phenomenologically \enquote{finite}. He calls Vop\v enka's infinite sets \enquote{inaccessible}. 
} 

It is possible to prove most of the usual properties of sets. %A linearly ordered set has the least and the greatest element. 
Lower-case letters denote sets from the set universe.
\medskip

\centerline{Axioms for classes}
\begin{enumerate}
\item[(B1)] \emph{Existence of classes:} If $\psi$ is \emph{any} formula then $\{x; \psi(x)\}$ is a class.
\item[(B2)] \emph{Extensionality:} Classes with the same elements are equal.
\item[(B3)] \emph{Existence of semisets:} There is a semiset. 
\end{enumerate} 

Classes are defined by the properties of sets from the universe of sets. \emph{Extended universe} is formed by these classes. They are denoted by capitals. 
\begin{definition} Let $X$ be a class from the extended universe. 
\begin{itemize} 

\item $X$ is a \emph{semiset} if it is not a set and it is a subclass of some set. 
\item $X$ is \emph{set-theoretically definable} if $X = \{x; \varphi(x)\}$ where $\varphi$ is a set formula 
\end{itemize}
\end{definition}

Each set $x$ is a set-theoretically definable class, since $x = \{y; y \in x\}$. The universal class $V = \{x; x = x\}$ is a set-theoretically definable proper class as well as Russell class $ \{x; x \notin x\}$. 

A set-theoretically definable subclass of a set is also a set.\footnote{Let $\varphi(x)$ be a set formula that defines the class $X = \{x; \varphi(x)\}$ and $X \subseteq a$. We wish to prove there is a set $b$ such that $b = X$. The set formula $(\exists b)(\forall x)(x \in b \Leftrightarrow (x \in a \wedge \varphi(x))$ is valid for $\emptyset$ and if it is valid for $a$ then it is valid for its successor $a \cup \{c\}$. According to the Induction axiom it is valid for all $a \in V$.} Hence, semisets cannot be set-theoretically definable.

Classes outside the extended universe can be coded within the extended universe. The problems concerning more general sets and classes are reduced to problems concerning the set universe and the extended universe.

For instance, all hairs of a hairy man can be coded by a set. The class of codes of hairs one can pull out so that the man is not bold is a proper subclass of a set of codes of all his hairs, a semiset from the extended universe. It is not set-theoretically definable; the Induction axiom does not lead to paradox.

\subsection{Finite and Infinite Classes}

The notion of the infinite is figuratively introduced in connection with a horizon. Its presence is represented by semisets. The infinite is comprehended as a phenomenon involved in the observation of large, incomprehensible classes. If we find a part that is a semiset in such a whole, we encounter the phenomenon of the infinite.

In this figurative sense, there are infinitely many people in the world, hairs of a hairy man, orange flowers, ancestors of a particular person, etc. 
\begin{definition}\label{fininf}
Let $X$ be a class. 
\begin{itemize}
\item $X$ is \emph{finite} if all its subclasses are sets. We denote it $Fin(X)$. 
\item $X$ is \emph{infinite} if it is not finite; some of its subclasses are a semiset, 
\end{itemize}

\end{definition}

This kind of infinity is different from the actual infinity in Cantor's set theory. The phenomenological meaning of infinity is primarily \enquote{the absence of easy survey}.

\subsection{Countable Classes}

The simplest phenomenon of infinite proper classes in the real world is described as a \enquote{path toward the horizon}. In every step of construction, we can make one further step; it is an unlimited process. 

A typical example is railroad ties that lead toward the horizon of our observation, the number of our male or female, ancestors or the ever-smaller reflections in two mirrors that face each other.

In AST, a \enquote{path toward the horizon} is formalized by the notion of a \emph{countable} class. A countable class is an infinite linearly ordered class such that all its initial segments are finite.

\begin{definition}
A class $X$ is \emph{countable}, if the following conditions are valid.
\begin{enumerate}
\item $X$ is infinite.
\item $X$ is linearly ordered by a relation $\leq$.
\item $(\forall x)(x \in X \Rightarrow Fin (\{y \in X; y \leq x\}).$
\end{enumerate}
\end{definition}

A countable class  does not have the greatest element; otherwise, it would be finite. It is neither set-theoretically definable; otherwise, it would have the greatest element. A linear ordering that defines a countable set is well-ordering.\footnote{Every subclass has the smallest element.} 
Thus any two countable classes are isomorphic with respect to the ordering.\footnote{It is easy to construct an isomorphism. We assign the smallest elements of both classes to each other. The remained subclasses have the smallest elements again, we assign them again to each other, and so on.}

\subsection{Prolongation Axiom}

People have always tried to go beyond the horizon; doing so is a typical human aspiration. 
A \enquote{path toward the horizon} does not stop exactly on the horizon. The observed object smoothly continues beyond the horizon. The illuminated part of an observed object can be extended at least a little more.

Railroad ties continue at least a bit beyond the horizon of observation. A boat on the sea continues its voyage beyond the horizon. A family tree has more male ancestors than we know about. This principle is incorporated into 

\begin{quote} (B4) \emph{Prolongation Axiom:} For each countable function $F$ there is a set function $f$ such that $F \subseteq f$.\end{quote}
It follows that every countable class $C$ can be extended to a set $c$ whose elements have the same set-theoretical properties, $C \subseteq c$. Hence, all countable classes are semisets because they are subclasses of infinite sets.

\section{Number Structures}\label{numb}

\subsection{Natural Numbers}

Natural numbers are constructed by the usual von Neumann way. Numbers are defined by a set formula and represented by the sets of smaller numbers from the universe of sets. Every number has a successor; there is no greatest natural number. Hence, the class of all natural numbers $\mathbb N$ is a set-theoretically definable infinite proper class.\footnote{It is not a set because the induction entails that a linearly ordered set-theoretically definable set has a least and the greatest element.} It is a model of Peano Arithmetic. 

Natural numbers that are finite in AST sense are called finite natural numbers, $\mathbb{FN}$. They lie before the horizon; they are \enquote{close, small, accessible}. If $n$ is finite, i.e. none of its subclasses is a semiset, then $n+1$ is finite too. So, there is no greatest finite natural number; the class $\mathbb{FN}$ is infinite, linearly ordered, and every initial segment is finite. Thus, $\mathbb{FN}$ is a countable class. It is a prominent mathematical representative for a \enquote{path toward the horizon}. 

Prolongation axiom guarantees the existence of an infinite natural number. Accordingly, there is a specific set $\alpha$ such that $\mathbb{FN} \subseteq \alpha$ and all its elements are natural numbers. Therefore, $\alpha$ itself is an infinite number that is greater than all finite numbers. The class $\mathbb{FN}$ is a proper cut on $\mathbb N$. 
$$ \mathbb{FN} = \{n \in \mathbb N; Fin(n)\} \subset \mathbb N.$$

The mathematical induction on $\mathbb N$ applies only to \enquote{definite} predicates that are characterised by set formulas. \enquote{Vague} predicates are characterised by any formula. (Dean 2018, p. 310). Such induction is valid only on $\mathbb{FN}$. 

This remark is a solution to paradoxes sorites. 
For example, Wang's paradox considers a vague predicate \enquote{small} that 
cannot be described by a set formula. Surely 0 is small. If $n$ is small then $n+1$ is small. Therefore, every number in front of the horizon (of being small) is small.

\subsection{Rational Numbers}

The \emph{rational} numbers $\mathbb Q$ are also constructed in the usual way as a quotient field over natural numbers $\mathbb N$. It is a set-theoretically definable, dense, linearly ordered field. Since $\mathbb N$ contain infinite numbers, for instance $\alpha \in \mathbb N \setminus \mathbb {FN}$, $\mathbb Q$ contain their inverses, ${1 \over \alpha}$, infinitely small numbers. 

\begin{definition}\label{is} Let $x \in \mathbb Q$ be a rational number. Then 
\begin{enumerate}
\item $x$ is \emph{infinitely small} if $(\forall n)(n \in \mathbb{FN} \Rightarrow |x| < {1 \over n})$; 
\item $x$ is \emph{infinite} if $(\forall n)(n \in \mathbb{FN} \Rightarrow |x| > n)$;
\item $x$ is \emph{bounded} if it is not infinite; the class of bounded numbers
$$ \mathbb {BQ} = \{x \in \mathbb Q; (\exists n \in {FN})|x| < n\};$$ 
\item $x, y \in \mathbb{BQ}$ are \emph{infinitely near} if their difference is infinitely small
$$x \doteq y \Leftrightarrow (\forall n)(n \in \mathbb{FN} \Rightarrow |x-y| < {1 \over n}).$$
\end{enumerate}
\end{definition}
Infinitely small numbers are beyond the horizon of a \enquote{depth}; infinite numbers are beyond the horizon of a \enquote{distance}; bounded numbers $\mathbb{BQ}$ lie in front of the horizon of a \enquote{distance}. 
Two bounded rational numbers are infinitely near if we do not distinguish their difference, resp. if we do not care about the difference between them.

The definition of infinitely small numbers assumes that the two horizons of a \enquote{depth} and of a \enquote{distance} correspond with each other. A number $n$ is infinite if and only if $\frac{1}{n}$ is infinitely small.\footnote{This agrees with Pascal's concept of two infinities: the infinitely large and the infinitely small. While they are infinitely distinct, they correspond to one another: from the knowledge of one follows the knowledge of the other. 
\enquote{La principale comprend les deux infinites qui se rencontrent dans toutes: l'une de grandeur, l'autre de petitesse. \dots Ces deux infinis, quoique infiniment diff\'{e}rents, sont n\'{e}anmoins relatifs l'un \` a l'autre, de telle sorte que la connaissance de l'un m\` {e}ne n\'{e}cessairement \` a la connaissance de l'autre. (Pascal 1866. pp. 288, 295).}} 

\subsection{Real Numbers}

Real numbers in AST express how we treat numbers in the life-world when they are associated with concrete objects. 
We talk about real numbers, but in most cases, we do not use them exactly. If we calculate the area of a real circle, we use a rational number that is as near to $\pi$ as we wish. 
If we divide a cake into eight equal pieces, each piece will probably be slightly larger or smaller than exactly one-eighth. When we speak about a trip twenty kilometres long, it is probably twenty kilometres plus or minus some meters.

A real magnitude commonly means a vague interval of rational numbers limited by a horizon. Real numbers are defined with the help of the relation of infinite nearness $\doteq$ on $\mathbb{BQ}$. A real number is represented by the class of infinitely near rational numbers that is called a monad.\footnote{The term \emph{monad} was originally borrowed by Robinson from Leibniz in his \emph{Nonstandard Analysis}. Vop\v enka took it from Robinson, and he had used it in the same meaning.}

\begin{definition}\label{ic}
A monad of $x \in \mathbb{BQ}$ is defined as a class of infinitely near rational numbers. 
$$Mon(x) = \{y \in \mathbb{BQ}; y \doteq x \}.$$
Real numbers $\mathbb R$ are represented by the class of monads. 
$$\mathbb R = \{Mon(x); x \in \mathbb{BQ}\}.$$
\end{definition}
The infinite nearness is reflexive and symmetric by Definition \ref{is}. In special cases, it is also transitive; it is an equivalence relation. By the factorisation of the infinite nearness, we obtain the class of real numbers.
$$\mathbb R = \mathbb{BQ}/ \doteq$$
Monads are factor-classes of this equivalence. They represent real numbers and have all the properties of classical real numbers. See discussion in Section \ref{non-standard}.

\section{Continuum in AST}\label{cont}

\subsection{$\sigma$-classes and $\pi$-classes}\label{sigma}

Before we describe the continuum phenomenon we shall define special classes characterised by \emph{phenomenal} predicates such as \enquote{looks red}, \enquote{tastes sour} or \enquote{sounds loud} (Dean 2018, p. 319). These phenomena are called \emph{primarily evident}.

Consider, for instance, \enquote{redness}. 
Some objects are definitely red, and some are less red; some are even
less red, and then there are borderline cases that balance between being red and being some
other colours: orange, violet or pink. 
The vagueness of redness is not something that can be somehow avoided. If we do not grasp redness in its vagueness, then we do not grasp it at all.
The vague boundaries of redness can be interpreted as the horizon before which the red colour is observed. 

Shades of redness can be precisely determined by
their wavelengths or pixels. The class of red objects can then be described
as a countable union of classes that are defined by specific wavelengths.

Similar points would seem to apply to other cases of primarily evident
phenomena. Although they are connected with vague properties, such phenomena can often be
described as the union of precisely measured objects. Primarily evident phenomena
will be described as \emph{$\sigma$-classes}. 

Then some classes are defined by the \emph{negations of primarily evident properties}.
Their existence is less conspicuous, and they often even do not
have their own names. They describe the borderline cases that connect primarily
evident properties. These phenomena will be described as \emph{$\pi$-classes}.

\begin{definition} A class $A$ is a \emph{$\sigma$-class} if it is a countable union of set-theoretically definable classes $A_n$. 
$$A = \bigcup\{A_n; n \in \mathbb{FN}\} \textrm{, \quad resp. } x \in A \Leftrightarrow (\exists n)(n \in \mathbb {FN}\wedge x \in A_n)$$
A class $B$ is a \emph{$\pi$-class} if it is a countable intersection of set-theoretically definable classes $B_n$. 
$$B = \bigcap\{B_n; n \in \mathbb{FN}\} \textrm{, \quad resp. } x \in B \Leftrightarrow (\forall n)(n \in \mathbb{FN} \impl x \in B_n)$$ 
\end{definition}

\enquote{It is cold} and \enquote{it is hot} are two distinct mutually exclusive and primarily evident phenomena. They can be represented by two disjoint $\sigma$-classes defined on the Celsius or Fahrenheit scale. Their complement is a $\pi$-class. Its interpretation expresses \enquote{it is fine}. This means \enquote{it is neither cold nor hot}. This phenomenon is not so striking, and often we are not completely aware of it. 

Countable classes, for instance $\mathbb{FN}$, are $\sigma$-classes.\footnote{If $C = \{c_1, c_2, c_3, \dots\}$ then $x \in C \Leftrightarrow (\exists n)(n \in \mathbb{FN} \wedge x \in C_n)$.} 
Both infinitely small and infinite numbers form $\pi$-classes.\footnote{
$x \in \mathbb Q$ is infinitely small iff $(\forall n)(n \in \mathbb{FN} \impl |x|< \frac{1}{n})$; great iff $(\forall n)(n \in \mathbb{FN} \impl |x|> n)$.}
Relations can be also $\sigma$-classes or $\pi$-classes depending on whether they are primarily evident or not. 
Infinite nearness $\doteq$ is a $\pi$-class.\footnote{If $x,y \in \mathbb{BQ}$ then $x \doteq y \Leftrightarrow (\forall n)(n \in \mathbb{FN} \impl |x - y| < \frac{1}{n})$.}

\subsection{Indiscernibility}\label{indiscernibility}

The original notion of continuum in antiquity, \emph{syneches}, meant anything around
us that we perceive as continuous. Space, movement and time are classic examples.
But we can also consider temperature, music, rainbow, etc., continuous. 

When treating continuum phenomena mathematically, AST accepts hypotheses 
that a continuum is produced when we observe a large but
remote class, and we are unable to distinguish its individual elements. They are indiscernible because they lie beyond the horizon of our observational capability.

For example, a heap of sand from a sufficient distance appears to be continuous. The indiscernibility of individual grains of sand is caused by the imperfection of our senses. However, this weakness is our advantage. Only due to the horizons of our perception do we perceive the world around us as continuous and cohesive. \enquote{The world without the horizon is unimaginable, maybe even impossible.} (Van Puersen 1977, p. 184).

While \emph{discernibility} of two objects is a primarily evident property - two objects are \emph{discernible} if there is an exact criterion that distinguishes them - and thus 
discernibility relation would be a $\sigma$-class; its negation \emph{indiscernibility} would be a $\pi$-class. Two objects are indiscernible if none of the exact criteria can distinguish them. Two grains of sand are indiscernible if their distance is less than any discernible distance.

\begin{definition} 
A binary relation $\approx$ defined on $V$ is an \emph{indiscernibility} if it is a countable intersection of set-theoretically definable, reflexive and symmetric relations $R_n \subseteq V \times V$ such that $R_{n+1}\subset R_n$. 
$$x \approx y \Leftrightarrow (\forall n)(n \in \mathbb{FN} \impl [x,y] \in R_n).$$
\end{definition}

By definition, indiscernibility is also a reflexive and symmetric relation and a $\pi$-class. A prominent mathematical example of indiscernibility is the infinite nearness $\doteq$ defined on bounded rational numbers $\mathbb{BQ}$.

\subsection{Continuum}
A continuum is understood as an infinite class endowed with an indiscernibility relation. For a continuous shape, it should be connected. 
\begin{definition}
Let $\approx$ be an indiscernibility relation defined on a set $w$. 
\begin{enumerate}
\item A class $X \subseteq V$ is \emph{connected} with respect to $\approx$ if for each non-empty proper subset $v \subset X$ there are indiscernible elements in $v$ and in its complement.
$$(\forall v)((v \subset X \wedge v \neq \emptyset) \impl (\exists x)(\exists y)(x \in v \wedge y \in X \setminus v \wedge x \approx y))$$ 
\item A \emph{monad} of $x \in V$ is the class of all elements indiscernible from $x$. 
$$Mon_\approx(x) = \{y \in V; x \approx y) \}$$
\item A \emph{figure} of $X \subseteq V$ is the class of all elements indiscernible from elements of $X$. 
$$Fig_\approx(X) = \{y \in V; (\exists x)(x \in X \wedge x \approx y) \}$$
\end{enumerate}
\end{definition}

A monad is a $\pi$-class. It lies precisely on the horizon of observation: it is a trace that $x$ left on the horizon. It contains indiscernible elements; it has no observable size. 

The figure of $X$ is a class of all indiscernible elements; it is a trace that $X$ left on the horizon, a shape of $X$ that we observe. It is simultaneously a union of all its monads. If $X$ is connected, then we perceive its figure as a continuous shape, a continuum.\footnote{This continuum concept can serve as a response to Zeno. Zeno's paradoxes are designed to refute both Aristotle's and Democritus'
views (Fletcher, p. 567). Among other things, he challenged the notion of the continuum as a plurality of things. He argued that if there are many things, then they need not have any size at all; otherwise, there would be unlimited objects. If things
have no size, then they do not exist at all.

Monads do not have observable size. But they are something:
They have a {body}. Joining or removing a monad is indistinguishable for an observer. However, the composition of infinitely many monads forms an observable part of the continuum.} 
$$Fig_\approx(X) = \bigcup\{Mon_\approx(x); x \in X\}$$

In some cases, the underlying class can be described directly. A heap of sand is
composed of grains. Two grains are indiscernible if their angular distance is less than circa one arcminute. 

Due to various measurement theories, almost all continua can be underpinned by
an infinite grid of coordinates, which is denoted by rational numbers dense enough to make very close coordinates indiscernible. The grid can be arbitrarily fine.

A ruler is applied to construct the sort of grid just described to a straight line drawn on a paper. It depends on our distance from the paper, which coordinates will appear indiscernible. 
A colour spectrum is measured by wavelengths or pixels. If the difference between two shades is less than \emph{circa} 4 nanometres and less, then it is indiscernible.

\subsubsection{Space}\label{space}
I will add a few of my expanding remarks on classic examples of a continuum: space, time and motion. The n-dimensional space can be described by a Cartesian system of coordinates denoted by n-tuples of rational numbers, $S \subseteq \mathbb Q^n$. Indiscernibility $\approx$ is a binary relation defined on $S$. In a concrete situation, it depends on our position, distance from an observable object and a chosen unit. Monads correspond to what is usually called \emph{points}. 

The simplest case is a surface from which we are always at the same distance $d$. Then the indiscernibility coincides with the infinite nearness, though it also depends on our distance from the surface.\footnote{\emph{Unit distance} means that the horizon of \enquote{depth} corresponds to the horizon of \enquote{distance}.} 
$$x \approx y \Leftrightarrow x/d \doteq y/d$$

Another case is when we stand in one place and look around. Two coordinates $x,y \in S$ are indiscernible if their difference is infinitely small due to their distance. 
$$x \approx y \Leftrightarrow \frac{x}{y} \doteq 1 \Leftrightarrow \frac{x - y }{y} \doteq 0$$
In this case, the size of the monads changes continuously. Nearby monads are very small; the more distant they are, the larger they are.

\subsubsection{Time}\label{time}

The choice of a unit of time coordinate depends on our interest; whether we deal with history, everyday life or maybe a race.
Time coordinates are represented by an infinite linearly
ordered subclass $T \subseteq \mathbb Q$. Two coordinates $x, y \in T$ are indiscernible, $x \sim y$, if the time interval between them is beyond the horizon of our perception. 
Monads of time are called \emph{instants}. 
For $x \in T$ $$Mon_\sim(x) = \{y \in T; y \sim x\}.$$

Since the relation $\sim$ is not generally transitive, instants are not disjoint, and they blend. We perceive time continuously as a flow of instants $Fig_ \sim(T)$. 
According to Vop\v enka, the \emph{past} and the \emph{future} are primarily evident phenomena. They can be described as $\sigma$-classes defined on $T$. The instant \emph{now} is a mere boundary between the past and the future, their complement; an instantaneous state that is a $\pi$-class.\footnote{\enquote{The past is that which has been present, the future that which will be present. So there cannot be either a past or a future unless there is, independently of past or future, such a thing as how things are now.} (Dummett 2000, p. 501).}

\subsubsection{Motion}\label{motion}\label{motion}

\emph{motion} is explained in AST as a phenomenon that we perceive when we are presented with a sequence of states in which the following state is indiscernibly from the preceding one. Like a movie that is a rapid series of images so that the differences between successive images lie beyond the horizon of our ability to distinguish. 

The simplest kind is a {motion of a point}. A point moves by changing permanently and imperceptibly its position. It is defined as a function $f$ from time $(T, \sim)$ to space $(S, \approx)$. 
$$f: (T, \sim) \longrightarrow (S, \approx)$$

A motion is \emph{continuous} if for any two indiscernible coordinates of time $x,y \in T$, their functional values are indiscernible too.
$$x \sim y \Rightarrow f(x) \approx f(y).\footnote{This definition corresponds to that of a uniformly continuous function in non-standard analysis. (Albeverio , p. 27)}$$

A motion is \emph{observable} if for any two discernible coordinates $x,y \in T$ their functional values are discernible.
$$ f(x) \approx f(y) \impl x \sim y.$$

What is commonly called a motion is both observable and continuous. The growth of a tree from a seed is also a continuous motion,
although it is not observable.

\section{Mathematical Idealization}\label{ideal}

The key question is, \enquote{Where is the horizon?} The provided examples are always joined with a specific real situation, such as grains, hairs or Celsius degrees. The horizon moreover depends
on the observer, their position and the chosen unit.
However, a horizon's mathematical form should be general and independent of
specific examples. And, of course, consistent.

\subsection{Abstraction and Idealization}
Vop\v enka devoted a part of his more philosophical book to the phenomenological description of the process of abstraction, although he did not call it that. He instead talks about \enquote{pulling new notions out of the maze of phenomena of the natural world}. 

\emph{Abstraction} can be described as incompleteness or a representation that lacks detail, but it is not designed to make literally false statements (Levy 2018, p. 7). All notions we introduced in Sections \ref{inf} and \ref{cont} of class, set, semiset, infinite set, countable class, $\sigma$-class and $\pi$-class, indiscernibility etc. arose as \emph{abstractions} of phenomena of the real world.

Abstraction concerns a description's degree of detail whereas \emph{idealisation} consists in
introducing simplifying misrepresentations. 
Idealisation affords to change some aspects of an object to obtain its ideal limit form. 
It is a \enquote{deliberate misrepresentation of some aspect of the world}. 
(Levy 2018, p. 7).\footnote{Levy gives an example: An account of gene flow in a population that assumes an infinite population size is idealised, in that; obviously, no real-world population is infinite.}

Garrison distinguishes two distinct directions in the process of idealisation: one ascending from the life-world, the other descending and applying to it. The first process of idealisation terminates its upward
movement with the objective, self-identical, universal forms, exact limit shapes. The second process descends \enquote{from the world of idealities to the
empirically intuited world.} (Garrison 1986, p. 330).

In the \emph{life-world}, a \enquote{path toward the hrizon} always
has its upper bound, i.e. a specific number that depends on concrete circumstances and that is infinite in AST sense but finite in the classical sense. 
Its \emph{abstract} form is the class $\mathbb{FN}$ of finite natural numbers. We have left open the question of whether this class is bounded by a specific number. 
However, the class of \emph{ideal} finite natural numbers has no such boundary and is closed under
arithmetic operations. It is \enquote{potentially infinite}. 
It is the exact limit shape, the objective, self-identical, universal form that we obtained by a deliberate misrepresentation

\subsection{The Non-Standard Model}\label{non-standard}

Ideal finite natural numbers $\mathbb{FN}$ are thus closed to addition and multiplication, and they 
form a model of Peano arithmetic. They can be represented as \emph{standard natural numbers} in a nonstandard model of PA. 
The prolongation axiom guarantees the existence of infinite natural numbers, which are represented as \emph{nonstandard natural numbers}. Formally, the mathematical model for finite and infinite natural numbers in AST can be considered an $\omega_1$-saturated nonstandard model of the natural numbers.\footnote{See for instance (Boolos, Burgess, and Jeffrey 2002. pp. 302 - 312).} This model guarantees the consistency of AST. 

Rational numbers $\mathbb Q$ form a dense, linearly ordered, non-Archimedean field. If $\mathbb{FN}$ form a model of Peano arithmetic, then the bounded rational numbers $\mathbb{BQ} \subseteq \mathbb Q$ form a commutative, linearly ordered ring. The relation of infinite nearness $\doteq$ 
defined on $\mathbb{BQ}$ is transitive.\footnote{
Let $x \doteq y$ and $y \doteq z$. Then for any $n \in \mathbb{FN}$ it is true that $|x - y| < {1 \over n}$. $\mathbb{FN}$ is closed under arithmetic operations, so also $|x - y| < {1 \over {2n}}$. The same holds true for $y \doteq z$. Consequently, $|x - z| \leq |x -y| + |y - z| < {1 \over {2n}} + {1 \over {2n}} = {1 \over n}$. Thus $x \doteq z$.} Hence, it 
is an equivalence, monads are disjoint equivalence-classes. 
By the factorization 
of $\mathbb{BQ}$ modulo $\doteq$ we obtain a linearly ordered field $\mathbb R$.\footnote{
The construction of \emph{real} numbers from non-standard \emph{rational} numbers is described in (Albeverio 1986, p.14).
Let $Q^*$ denotes the class of nonstandard rational numbers, 
it is a dense, linearly ordered, non-Archimedean field. 
Let $\mathbb Q_b$ denote the set of bounded rational numbers, $\mathbb Q_i$ the set of infinitely small rational numbers, 
$\mathbb Q_i \subseteq \mathbb Q_b \subseteq \mathbb Q^*.$ $\mathbb Q_i$ form the maximal ideal in a ring $\mathbb Q_b$. 
The result of the factorisation of $\mathbb Q_b$ modulo $\mathbb Q_i$ is the same as a factorisation by the infinite nearness $\doteq$. We obtain the field isomorphic to real numbers. 
$$\mathbb Q_b /\mathbb Q_i \quad = \quad \mathbb Q_b/\doteq \quad \cong \quad \mathbb R.$$
The class $Q^*$ corresponds to $\mathbb Q$ of AST, $\mathbb Q_b$ to $\mathbb{BQ}$, infinite nearness has the same definition.}
$$\mathbb{BQ}/\doteq \quad \cong \quad \mathbb R$$
The class $\mathbb R$ has all the properties of usual real numbers. 
It is a linearly ordered, dense, complete, and Archimedean field.

\subsection{Limits of Idealization}

However, the horizon in the real world is always classically finite.\footnote{Discussion of a similar systems can be found in (Dean 2018, pp. 309 - 313). Dean inquires models of a theory $S_\tau$ formulated over a language extending that of first-order arithmetic with a new predicate $F(x)$ such that $F(0) \wedge (F(x) \impl F(S(x)) \wedge (F(x) \impl (\forall y)(y < x \impl F(y))$ but $\neg F(\tau)$ for a sufficiently great term $\tau$. An interpretation of $F(x)$ is $x$ is \emph{feasible}. It can express any soritical predicate. If we interpret it as $x$ \emph{is finite} then $S_\tau$ and $\mathbb{FN}$ have the same models.
Dean suggested the \emph{neo-feasibilist theory of vagueness} as a possible solution. This theory employs a nonstandard model of natural numbers. The term $\tau$, which represents a non-feasible number, is realised by an infinite number, and the soritical predicate is interpreted as a proper cut on natural numbers). Dean cited Vop\v{e}nka as the only person to use nonstandard methods in connection with vagueness. (Dean 2018, p. 296).}
Consequently, the indiscernibility commonly is not transitive, which matches our experience. In the same context,
Dummett speaks about the \enquote{non-transitivity of the relation of not discriminable difference} as one feature of reality - or of our experience with it.
Take, for example, the rainbow. We can always divide it into narrow strips so that two neighbouring strips are indistinguishable or appear to be of the same colour. However, after adding a finite number of strips, we can see that the colour has changed. It is similar to other continua, for instance, space or motion.\footnote{\enquote{I look at something which is moving, but moving too slowly for me to be able to see that it is moving. After one second, it still looks to me as though it was in the same position; similarly, after three seconds. After four seconds, however, I can recognise that it has moved from where it was at the start, i.e. four seconds ago.} (Dummett 315).} 
This implies that the monads of the observed continuum are not generally disjunctive but overlap. We do not distinguish them individually, although this is precisely why we perceive the structure as continuous.
\footnote{The most famous Zeno's paradoxes are based on the tension between the real continuum and the ideal mathematical continuum. 
According to \emph{Dichotomy paradox} 
\enquote{there is no motion because that which is moving must reach the midpoint before the end.} (McKirahan, p. 181). 
Since the argument can be repeated again and again, one must go through infinitely many places before arriving at the goal. No finite distance can ever be travelled: all motion is impossible. 
Indeed, an interval of ideal real numbers that has the length one can be halved
again and again, and still, it is an interval of real numbers. The distance from
the end will subsequently be
$\frac{1}{2},\frac{1}{2^2}, \frac{1}{2^3}, \dots, \frac{1}{2^n}, \frac{1}{2^{n + 1}}, \dots$
but never $0$. It is impossible to reach the goal.
However, in a concrete situation, there is a finite number $n$ such that the distance $ \frac{1}{2^n}$ is indiscernible from the endpoint, $\frac{1}{2^n} \approx 0$. We are at the goal, at the same monad, in $n$ steps
The \emph{Achilles and the tortoise} paradox is based on the same principle. Their distance becomes indiscernible after finitely many steps. They are within the same monad. If the race continues, Achilles and the tortoise go on
from the same position.}

What can we do? First, we can investigate the ideal mathematical description,
i.e. a model of a nonstandard theory, and carefully apply its results to the real world. 
The boundaries of idealisation are known. 
\begin{quote}The idealities which are assumed by the descending process of idealisation serve as a priori \enquote{guides} in the further determination of those vague and inexact empirical entities within the manifold of intuition. \dots
The vague 
entities of everyday experience are
geometrically determinable to precisely the degree that they \enquote{participate}
in the pure geometrical forms which they motivate. (Garrison 1986, p. 331) \end{quote}

To a large extent, this is exactly what infinitary mathematics does and what owes its success even without being aware of it. 
\begin{quote} 
We apply the mathematical results in the real world -- in the same way the results of ancient geometry have been applied ever since its origin. That is by replacing the ideal horizon with the horizon limiting the human sight into the real world. In doing so, we evidence the inevitable misrepresentations that it brings.
(Vop\v enka 2015, p. 89). \end{quote}

Second, we can investigate the abstract concepts without their idealisation. Many mathematical statements can be proclaimed on $\pi$-classes, $\sigma$-classes, indiscernibility, monads etc., as we introduced them in Sections \ref{inf} and \ref{cont}. 
We can carefully weaken traditional mathematical concepts and replace them with more convenient ones. It opens a broad field of applications.   %That is what Vop\v enka tried to do in his last work NIM. 

Third, the phenomenological interpretation opens possibilities for new concepts. It enables 
a mathematical treatment of notions that either have not yet been defined mathematically or that have been defined otherwise. 
New questions open up. We can ask about precise conditions of continuous motion, different space indiscernibilities or a horizon shift. Sections \ref{space}, \ref{time} and \ref{motion} are an example of the direction that could be pursued.

\subsection{The Limit Universe and the Witnessed Universe}

Vop\v enka distinguished the study of a \emph{witnessed universe} - when we guarantee the existence of a semiset included in a specific set - from the study of a \emph{limit universe} - when subsemisets of all sets are eliminated. 
The witnessed universe is \emph{abstracted} from the real world, whereas the limit universe is its \emph{idealisation}.

In (Vop\v{e}nka 1979), he decided to study the limit universe because 
\enquote{the problems of the witnessed universe were not yet satisfactorily understood.} Apparently, he meant Yesenin-Volpin's feasible numbers because studying them is akin to studying the witnessed universe. But he admitted that \enquote{all situations to which our theory applies in the world perceivable by our senses correspond to the witnessed universe}. That \enquote{such situations motivate various notions introduced in AST}. (Vop\v{e}nka 1979, p. 38).

However, this probably was the main reason why he repeatedly returned to his theory. In NIM, he attempted to investigate the proper cuts on natural numbers in the witnessed and the limit universe simultaneously.

\section{Conclusion}\label{concl}

Instead of searching for infinity
somewhere in the immense cosmic distances, AST places it into this world
where it becomes a figurative designation for large, incomprehensible collections.

If we observe a vague part in such a whole, then we encounter a phenomenon of
infinity. Vague parts are limited by the horizon. What is before the horizon is
finite; what is beyond the horizon is infinite.

The interpretation of continua is based on the
hypothesis that a continuum phenomenon is produced due to an observation
of an infinite but remote class such that its elements are indiscernible. 
We can liken it to an impressionist painting. It consists of an infinite
number of brushstrokes that are visible from close range but merge into one image from afar. Only because of the 
imperfection of our senses described by the indiscernibility relation 
we perceive the given class as continuous.

To capture this interpretation mathematically and express its \emph{abstract}
form, Vopenka introduced new notions: semiset, countable class, $\sigma$-class, $\pi$--class, indiscernibility, and so on. 

The meaning of the terms rational and real numbers is literal. While rational numbers are rationally constructed as a quotient field over natural numbers, real numbers express how we perceive the size of objects in specific real-world situations. They are defined as monads of indiscernible rational numbers. 

It is necessary to take one more step: \emph{idealisation}. %Of course, the horizon, as well as the indiscernibility in the life-world, depends on the observer, their position and the chosen unit. However, its mathematical form should be independent of a concrete situation. 
Finite natural numbers that represent the ideal \enquote{path toward a horizon} are unbounded and closed under arithmetic operations.
This implies that an ideal indiscernibility is transitive, ergo it is an equivalence, and its monads are disjoint. Specifically, monads of an infinite nearness defined on the bounded rational numbers 
have all the usual properties of real numbers.

$$\ast \ \ \ast \ \ \ast$$

Vop\v enka's theory seems to be closed to the feasibility theory. Indeed, finite numbers
in his witnessed universe resemble feasible numbers. Bellotti analyses the possibility of modelling the infinite in the finite. 
(Bellotii 2008. p. 3). According to him, there are only two reasonable ways to model countable infinity in the finite: either feasibility theory, which is only almost consistent, or nonstandard
methods. He concludes that although many interesting results have been obtained in these attempts, they ultimately show that no satisfactory reduction is possible.\footnote{\enquote{We have seen that when we consider the two main alternatives which apparently allow one to make sense of a sort of 'modelling' of countable infinity in the finite, namely nonstandard methods and feasibility, we face a dilemma. If we take into account proofs of arbitrary finite length, we might have consistency proofs, \dots, but we do not obtain any reduction of the infinite to the finite. On the other hand, if we consider proofs of length at most $k$, with $k$ a standard integer, we have only proofs (possibly in relatively weak theories) of 'almost consistency', and we do not obtain real consistency proofs.} (Bellotti 2008, p. 23).}

However, AST endeavour is a countermovement. 
He \enquote{models} finiteness by the infinite. Infinity is the idealisation of a great finite. 
One's handling ideal infinite objects is often more straightforward than one's handling real complex things. It also explains why results of infinite mathematics are applicable in the real world.

Vop\v enka derived his new notions from real-world phenomena.
He stated their abstract forms. Mathematical results can be achieved 
at this level, but other, deeper results can be obtained by the idealisation of a
horizon. AST demonstrates the process of idealisation and its
limits simultaneously.

The phenomenological interpretation of infinity and its connection with vagueness opens a broad field for applications of mathematical research. Phenomena of the real world that have not yet been possible to capture mathematically, such as motion, time, primarily evident properties, are described.   
It is no coincidence that AST also provides a natural solution
to some ancient puzzles such as sorites, paradoxes of motion and the composition
of a continuum. 

\bigskip

\textbf{Acknowledgements:} 
I am grateful to Alena Vencovsk\' a for reading the text and helping me with the mathematical questions. I have also used her prepared English translation of Vop\v eka's latest book \textit{New Infinitary Mathematics}. For valuable comments I thank Miroslav Hole\v cek and Pavel Zlato\v s. I would also like to thank three anonymous reviewers for all their comments and for their patience. They considerably helped to improve the text.

\end{document}